\documentclass[12pt,a4paper]{article}
\usepackage{amssymb}

\usepackage{graphicx,amssymb,amsfonts,epsfig,amsthm,a4,amsmath,url}
\usepackage[latin1]{inputenc}

\newtheorem*{thmm}{Theorem}

\newtheorem{thm}{Theorem}[section]

\theoremstyle{definition}

\theoremstyle{remark}
\newtheorem{rem}[thm]{Remark}

\numberwithin{equation}{section}

\newcommand{\N}{\mathbf{N}}

\newcommand{\bpr}{\noindent \textbf{Proof}: ~}

\newcommand{\epr}{~$\blacksquare$}

\title{The inclusion of the Schur algebra in $B(\ell^2)$ is not inverse-closed.}
\author{Romain Tessera\footnote{The author is partially supported by the
NSF grant DMS-0706486.} }
\date{\today}

\begin{document}

\baselineskip=16pt

\maketitle

\begin{abstract}
The Schur algebra is the algebra of operators which are bounded on $\ell^1$ and on $\ell^{\infty}$.  In \cite{Sun}, Sun conjectured that the Schur algebra is inverse-closed. In this note, we disprove this conjecture.
Precisely, we exhibit an operator in the Schur algebra, invertible in $\ell^2$, whose inverse is not bounded on $\ell^{1}$ nor on $\ell^{\infty}$. 
\end{abstract}

The Schur algebra is the unital algebra of infinite matrices whose rows and columns are uniformly bounded in $\ell^1$. Such matrices define operators which are uniformly bounded on $\ell^p$ for all $1\leq p\leq \infty$. In this short note, we prove the following

\begin{thmm}
There exists an infinite symmetric matrix $M=\{m_{i,j}\}_{i,j\in \N}$ such that 
\begin{itemize}

\item $m_{ij}=0$ or $1/4$, \item the support of each row and each column has cardinality $4$, \item $I-M$ is invertible in $\ell^2$, but not in $\ell^{\infty}$.
\end{itemize}
\end{thmm}
\bpr Let us consider a finitely generated group $G$, equipped with a probability measure $\mu$ on $G$ such that $\mu(g)=\mu(g^{-1})$ for all $g\in G$, and such that the support of $\mu$ is finite and generates the group $G$. Let $M$ be the operator of convolution by $\mu$ on $\ell^2(G)$, i.e.
$$M(f)(g)=\mu\ast f(g)= \sum_{h\in G} m(g^{-1}h)f(h)=\sum_{h\in G} m(h)f(gh).$$ Up to enumerate the elements of $G$, one can see $M$ as an infinite matrix.
Note that the cardinality of the support of both the rows and the columns of $M$ is simply the cardinality of the support of $\mu$.

Let us check that if $G$ is infinite, then the convolution by $\mu$ is never invertible in $\ell^{\infty}.$ Note that since it is self-adjoint, it is neither invertible on $\ell^1$. 
Let $S$ be the support of $\mu$. The {\it word metric} on $G$ is defined as follows: $d_S(g,h)=\inf\{n\in \N;  g^{-1}h=s_1\ldots s_n, s_i\in S\}$. The ball  $B(e,n)$, of radius $n$ and centered on the neutral element $e$ is therefore the set of all $g$ which can be written as a product of at most $n$ elements of $S$. 

For each $n\geq 1$, let  $f_n$ be the function measuring the distance to the complement of $B(e,n)$ in $G$, i.e. $$f_n(g)=\min_{h\in G\setminus B(e,n)}d(g,h).$$ 
Obviously, $f_n$ is a $1$-Lipschitz function on $(G,d_S)$. Therefore, by definition of $M$, one has that
$|(I-M)(f_n)(g)|\leq 1$ for all $g\in G$.
But on the other hand, $f_n(e)=n$, so we obtain the following inequality
$$\frac{\|(I-M)(f_n)\|_{\infty}}{\|f_n\|_{\infty}}\leq 1/n,$$
which tends to $0$ when $n\to \infty.$
Hence $I-M$ is not (left) invertible in $\ell^{\infty}.$

On the other hand, by a classical result of Kesten\footnote{Kesten gives a probabilistic proof of his result. For a more analytic approach, one can consult for instance \cite{C}.} \cite{Kest}, the group $G$ is non-amenable if and only $I-M$ is invertible in $\ell^2(G)$. 
The most classical example of a non-amenable group is the free group with two generators $\langle x,y\rangle$. Taking $\mu$ such that $\mu(x)=\mu(y)=1/4$, one gets the precise statement of the theorem.
\epr

\begin{rem}
In the case of the free group with $2$ generators, and for $\mu$ as above, one has $\|(I-M)^{-1}\|_2=2/(2-\sqrt{3})$ \cite{Kest}. Note that this norm is just the inverse of the smallest real eigenvalue of the so-called {\it discrete Laplacian} $\Delta= I-M$ on the Cayley graph of the free group $\langle x,y\rangle$ (which is a 4-regular tree). Another formulation of Kesten's theorem is that a finitely generated group is non-amenable if and only if any (resp. one) of its Cayley graphs has a spetral gap (i.e. this first eigenvalue is non-zero). 
\end{rem}

\bigskip
\footnotesize


\begin{thebibliography}{KM98b}

\bibitem[Cou]{C} T. {\sc Coulhon}.
\newblock {\em Random walks and geometry on infinite graphs.}
\newblock  Lecture notes on analysis on metric spaces,
Luigi Ambrosio, Francesco Serra Cassano, eds., 5-30, 2000.

\bibitem[Kest]{Kest} H. {\sc Kesten}. \newblock {\em
Symmetric random walks on groups.} \newblock Trans. Amer. Math. Soc.
92, 336-354, 1959.


\bibitem[ShS]{Sun} C. E. {\sc Shin} and Q. {\sc Sun}. \newblock  {\em Stability of localized operators}.  \newblock Journal of Functional Analysis,  256, 2417--2439, 2009.
%\bibitem[Ni]{Ni} B. {\sc Nica}. \newblock {\em A study of decay.}
%\newblock
%Preprint, 2007.

\end{thebibliography}
\end{document}